\documentclass[12pt,thmsa]{article}
\usepackage{amsmath}
\usepackage{amsfonts}
\usepackage{amssymb}

\usepackage[applemac]{inputenc}  

 \usepackage{amsmath, amsthm, amssymb, amsfonts, enumerate}

\def \R{I\!\!R}

\newtheorem{thm}{Theorem}
\newtheorem{theorem}[thm]{Theorem}

\newtheorem{lem}[thm]{Lemma}
\newtheorem{proposition}[thm]{Proposition}

\def \lim{ {\rm lim} }

\def \supp{ {\rm supp}\;  }
\def \argmax{ {\rm Argmax} }

\begin{document}

\title{Strategic information transmission with sender's approval}
\author{Fran\c{c}oise Forges\thanks{%
PSL, Universit\'{e} Paris-Dauphine, LEDa. E-mail: francoise.forges@gmail.com}
\ and J\'{e}r\^{o}me Renault\thanks{%
TSE, Universit\'{e} de Toulouse 1, D\'{e}partement de Math\'{e}matiques.}}
\date{January 2020\thanks{%
This research started during the winter 2015-2016, while the first author
was visiting Humboldt University, Berlin.\ We did not know of Shimizu (2013,
2017) at the time.\ These papers were pointed out to the first author by
Daniel Kr\"{a}hmer on the occasion of a talk in Bonn in April 2018. Then we
discovered Matthews (1989), among the references of Shimizu (2013, 2017).} 
\thanks{%
We thank Anna Bogomolnaia, Ulrich Horst, Vincent Iehl\'{e}, Fr\'{e}d\'{e}ric
Koessler, Daniel Kr\"{a}hmer, Ehud Lehrer, Ronny Razin, Antoine Salomon,
Roland Strausz and Bertrand Villeneuve for stimulating conversations on the
topic of this research. We also thank the participants of the workshops
\textquotedblleft Competition and Incentives\textquotedblright\ (Humboldt
University, Berlin, June 9-10, 2016), \textquotedblleft Game theory: \
Economics and Mathematics (GEM)\textquotedblright\ (South University of
Denmark at Odense, October 5-6, 2018), \textquotedblleft Signaling in
Markets, Auctions and Games: a Multidisciplinary Approach\textquotedblright\
(Paris 2, May 22-23, 2019), \textquotedblleft Frontiers in
Design\textquotedblright\ (University College London (UCL), June 14-15,
2019), \textquotedblleft Information Design and Splitting
Games\textquotedblright\ (Centre d'Economie de la Sorbonne, Paris, June
17-19, 2019), as well as the audiences of seminars in Besan\c{c}on, Bonn,
Lancaster, London (LSE), Paris-Dauphine and Rome (LUISS). J\'{e}r\^{o}me Renault
gratefully acknowledges funding from ANR-3IA Artificial and Natural Intelligence Toulouse Institute, 
grant ANR-17-EUR-0010 (Investissements d'Avenir program) and ANR MaSDOL.}}
\maketitle

\begin{center}
\newpage

\textbf{Abstract}
\end{center}

We consider a sender-receiver game with an outside option for the sender.\
After the cheap talk phase, the receiver makes a proposal to the sender,
which the latter can reject.\ We study situations in which the sender's
approval is crucial to the receiver.

We show that a partitional, (perfect Bayesian Nash) equilibrium exists if
the sender has only two types or if the receiver's preferences over
decisions do not depend on the type of the sender as long as the latter
participates. The result does not extend: we construct a counter-example
(with three types for the sender and type-dependent affine utility
functions) in which there is no mixed equilibrium. In the three type case,
we provide a full characterization of (possibly mediated) equilibria.\newpage

\section{Introduction}

We consider a general model of sender-receiver games.\ The specific feature
of our games is that the sender has an outside option. After the cheap talk
phase, the receiver proposes a decision to the sender; if the sender
approves it, the decision is made; otherwise, the sender chooses his outside
option, which can be interpreted as \textquotedblleft exit\textquotedblright
. Under complete information, the game reduces to an ultimatum game, in
which one player makes a \textquotedblleft take it or leave
it\textquotedblright\ offer to the other.\ In our framework, this other
player has private information and can send a costless message to the
receiver before getting an offer.

We are interested in situations in which the sender's approval is crucial to
the receiver.\ We thus assume that the receiver's utility in case of exit is
very low, as compared to what he can expect if the sender accepts his
proposal. It is not difficult to find examples in which a decision-maker
consults with an informed party before making a proposal that can be
ultimately rejected and in which rejection has unvaluable, damaging
consequences for the decision-maker. For instance, firms try to figure out
workers' requirements in order to avoid strikes and boycotts.\ Governments
discuss with kidnappers, hoping that hostages will not be killed. As a third
example, analyzed in Matthews (1989), the U.S.\ Congress may worry about the
President's veto.

We focus on equilibria in which the sender does not make uncredible threats
at the approval stage, namely, accepts a proposal if and only if it gives
him at least the utility of his outside option.\footnote{%
Our solution concept is basically subgame perfect Nash equilibrium. Except
for the approval stage, our model behaves as a standard cheap talk game, in
which Perfect Bayesian equilibrium is not restrictive.} We ask whether our
sender-receiver game has an equilibrium in which exit does not occur.

To answer this question, we introduce an auxiliary \textquotedblleft limit
game\textquotedblright\ $\Gamma $, in which equilibria are necessarily
without exit. The equilibria of $\Gamma $ are easily characterized by two
sets of conditions: incentive compatibility and constrained optimization.
Both sets of conditions are tractable but satisfying them jointly is
demanding. Existence of an equilibrium in the game $\Gamma $ is not
obvious.\ For instance, as opposed to standard sender-receiver games, $%
\Gamma $ may not have any nonrevealing equilibrium. This means that, in
absence of information transmission, the receiver cannot make any decision
that would give at least his reservation utility to the sender, whatever his
type. We identify various assumptions which guarantee that, in a situation
like this, the sender can credibly reveal some information to the receiver,
in such a way that exit will never happen.

We maintain the following assumptions on the game $\Gamma $: the sender has
finitely many types (which can be multidimensional, e.g., belong to $\mathbb{%
R}^{n_{1}}$, for some $n_{1}$), the receiver has a compact set of decisions
(typically, a closed, bounded set in $\mathbb{R}^{n_{2}}$, for some $n_{2}$)
and both players' utility functions are continuous.\footnote{%
This covers the particular case where the receiver has finitely many
actions, over which he can randomize.} We also make the \textit{sine qua non}
assumption that under complete information, i.e., when the receiver knows
the sender's type, there exists a decision that gives the sender at least
his reservation utility. We then consider the subsets of all types such that
there is a decision inducing them to participate and we call
\textquotedblleft participation structure\textquotedblright\ the maximal
subsets (with respect to set inclusion).\ For instance, if the sender has
only two possible types, 1 and 2, the participation structure is either $%
\left\{ \left\{ 1\right\} ,\left\{ 2\right\} \right\} $ or $\left\{
1,2\right\} $.

We establish that the game $\Gamma $ has a \textit{partitional} equilibrium,
namely, an equilibrium in which the sender's strategy is pure, in the
following cases:

\begin{description}
\item[(i)] the sender has two types;

\item[(ii)] the participation structure is a partition of the type set;

\item[(iii)] the decision set is a real interval and for every type, the
sender's utility function is monotonic in the receiver's decision;

\item[(iv)] the receiver's utility function -- when the sender participates
-- does not depend on the sender's type.
\end{description}

Cases (i) and (ii) are rather straightforward, with (ii) generalizing (i).
Case (iii) applies in particular when the receiver has only two actions,
over which he can randomize.\ Case (iv) has the most important scope.\ It
applies as soon as the receiver knows his own preferences over decisions,
but is eager to make a choice that will ensure the -- type-dependent --
informed player's participation. Existence of a partitional equilibrium
under (iv) is the main result of the paper (Theorem 8).

The previous assumptions may look restrictive, but, without them, existence
of an equilibrium in $\Gamma $ cannot be guaranteed, even if the sender is
allowed to use a mixed strategy.\ We indeed propose an example, in which the
sender has three types, the receiver has three actions, the participation
structure is not a partition and the utility functions are type-dependent.
In this example, there is no \textit{mixed} equilibrium.\ However, an
equilibrium does exist if the information transmission stage is handled by a
mediator.

Finally, we propose a complete analysis when the sender has three possible
types. We identify two kinds of \textquotedblleft participation
structures\textquotedblright\ beyond the straightforward case (ii) above.
The first one arises in the example mentioned in the previous paragraph.\
Existence of a mediated equilibrium can then be established.\ In the other
case, another example shows that there may not be any partitional
equilibrium.\ However we prove that there always exists a mixed equilibrium
if the decision set is convex and the utility functions are affine
(Proposition 9).

Here is a description of the paper.\ We discuss the related literature
below. In Section 2, we make the sender-receiver game $\Gamma $ and the
solution concept fully precise.\ Propositions 3 and 4 (which are established
in Section 6.1) allow us to argue that the game $\Gamma $ is relevant to our
study. In Section 3, we establish existence of a partitional equilibrium in $%
\Gamma $ under assumptions (i), (ii), (iii) or (iv) above. Our main result,
Theorem 8, is associated with case (iv).\ Section 4 is devoted to examples.\
Sections 4.1 and 4.2 propose a family of kidnapping games.\ Section 4.1
illustrates partitional equilibrium.\ Section 4.2 proposes a game that does
not have any partitional equilibrium but has a mixed equilibrium. Section
4.3 goes on with a game that does not have any mixed equilibrium but has a
mediated equilibrium. The three type case, including Proposition 9, is the
topic of Section 5. Section 6 is an appendix containing the proofs of
Propositions 3, 4 (Section 6.1) and 9 (Section 6.2).

\bigskip

\noindent \textbf{Related papers}

Shimizu (2013, 2017) adds an approval stage to Crawford and Sobel (1982)'s
sender-receiver game, in the popular case where the prior is uniform over
the unit interval and the utility functions are quadratic.\ He assumes, as
we do, that exit is damaging for the receiver but the setup is otherwise
quite different from ours. He shows that, in his particular model, credible
exit possibilities can make cheap talk informative even when the players'
conflict of interest is relatively large.

Matthews (1989) studies a sender-receiver game motivated by a specific
application, in which the sender is the U.S. President, the receiver is the
Congress and the decision is about a practical matter, like the level of
military expenditures. The President can veto the Congress' proposal.
Preferences are unimodal, as in Shimizu (2013, 2017), but the receiver's
utility does not depend on the sender's type (as in the current paper,
Theorem 8, Section 3.4).\ More importantly, in Matthews (1989)'s model, the
sender's rejection leads to status quo, rather than to exit, and does not
necessarily yield a very low utility to the receiver. Matthews (1989)'s
point is to show that thanks to incomplete information on the President's
type, veto can happen at equilibrium, i.e., without relying on uncredible
threats.

Our model can be viewed as a principal-agent problem in which the principal
-- alias the receiver -- cannot commit to a mechanism at the ex ante stage.\
This is an extreme case of Bester and Strausz (2001)'s principal-agent
problem with limited commitment.\ In this context, it makes sense to allow
the agent -- alias the sender -- to veto the principal's decision.\ Under a
mechanism design perspective, the principal looks for an equilibrium that
gives him the best ex ante expected utility, which amounts to an equilibrium
in which all types accept the principal's proposal if the principal's
utility, when the agent chooses his outside option, is sufficiently low.
This leads us to impose individual rationality conditions for the agent at
the \textquotedblleft posterior\textquotedblright\ stage, i.e., after the
principal has made a proposal. The relevance of posterior individual
rationality and its impact on incentive compatibility have been stressed in
a number of papers, e.g., Gresik (1991), Compte and Jehiel (2007, 2009)
Forges (1990, 1999) and Matthews and Postlewaite (1989).

Finally, Forges and Horst (2018)'s concept \textquotedblleft talk and
cooperate (perfect Bayesian) equilibrium\textquotedblright\ (TCE, Section
5.3) is motivated by the same questions as the present paper, but is defined
in a different model: the sender also has to make a decision, which is
relevant to his own payoff only.\ At a TCE, the receiver (who can be
interpreted as a principal) proposes a joint decision, which the sender
accepts whatever his type.\ Should player 1 reject player 2's proposal, both
players would choose an action, independently of each other. By contrast, in
the present paper, the sender just chooses an outside option.\ Forges and
Horst (2018) establish an existence result for another solution concept --
\textquotedblleft cooperate and talk (perfect Bayesian)
equilibrium\textquotedblright\ (CTE) -- but just indicate that their
methodology does not apply to TCE.

\section{Model}

\subsection{Sender-receiver games}

We start with a family of games\textbf{\ }$\Gamma (v_{0})$, $v_{0}\in 
\mathbb{R}$, between a sender (player 1) and a receiver (player 2).\ $\Gamma
(v_{0})$ is described as follows:

\begin{itemize}
\item A type $k\in K$ is chosen according to a prior probability $p\in
\Delta (K)$.

\item Player 1 is informed of $k$.

\item Player 1 sends a message $m\in M$ to player 2.

\item Player 2 proposes a decision $x\in X$ to player 1.

\item If player 1 accepts player 2's proposal, the decision $x$ is enforced,
player 1 gets $U^{k}(x)$ and player 2 gets $V^{k}(x)$.

\item If player 1 rejects player 2's proposal, player 1 chooses an outside
option and gets $u_{0}^{k}$.\ Player 2 gets $v_{0}$.
\end{itemize}

\bigskip

We assume that:

\begin{itemize}
\item The set of types $K$ is finite\footnote{%
We do not make any assumption beyond the fact that there are finitely many
types; in particular, types can be \textquotedblleft
multidimensional,\textquotedblright\ with $K\subset \mathbb{R}^{n_{1}}$, for
some $n_{1}\geq 1$.} and $p^{k}>0$ $\forall k\in K$.

\item The set of messages $M$ is finite, such that $\mid M\mid \geq \mid
K\mid $.

\item The set of decisions $X$ is compact.\footnote{%
We will indicate explicitly when $X$ will be required to be convex.} As a
typical example, $X\subset \mathbb{R}^{n_{2}}$, for some $n_{2}\geq 1$; for
instance, player 2 has a finite set of actions $A$ and $X=\Delta (A)$
corresponds to the set of mixed strategies of player 2.

\item The utility functions $U^{k}:X\rightarrow \mathbb{R}$ and $%
V^{k}:X\rightarrow \mathbb{R}$ are continuous; for instance, if $X=\Delta
(A) $, $U^{k}$ and $V^{k}$ may correspond to expected utility.
\end{itemize}

We further assume that:

\begin{itemize}
\item For every $k\in K$, there exists$\ x\in X\ $such that$\ U^{k}(x)\geq
u_{0}^{k}$.

\item For every $k\in K$, for every $x\in X$, $V^{k}(x)\geq v_{0}$, namely,

\begin{description}
\item $v_{0}\leq \min_{k\in K}\min_{x\in X}V^{k}(x)$.
\end{description}
\end{itemize}

\bigskip

We are interested in situations in which the sender's approval is crucial to
the receiver, namely, in which $v_{0}$ can be arbitrarily low.\ Let us
denote as $\Gamma $ the \textquotedblleft limit game,\textquotedblright\ in
which $v_{0}=-\infty $. We will show that $\Gamma $ is a tractable tool,
which is appropriate to study $\Gamma (v_{0})$ when $v_{0}$ is small enough.

\bigskip

Let us set, for every $L\subseteq K$%
\begin{equation}
X(L)=\left\{ x\in X:U^{k}(x)\geq u_{0}^{k}\text{, }k\in L\right\} \text{.}
\label{X(L)}
\end{equation}%
Given a subset of types $L$, $X(L)$ is the set of decisions that are
acceptable by all types in $L$. We write $X(k)$ for $X(\left\{ k\right\} )$,
so that $X(L)=\bigcap\limits_{k\in L}X(k)$.

\subsection{Equilibria}

Our solution concept, in $\Gamma (v_{0})$ and $\Gamma $, is basically
subgame perfect Nash equilibrium, but perfect Bayesian equilibrium would not
be more demanding: as in standard sender-receiver games, finding beliefs
rationalizing player 2's choices is not an issue. What is crucial here is to
avoid non-credible threats from player 1. In the sequel, we simply refer to
\textquotedblleft equilibrium.\textquotedblright

At a subgame perfect Nash equilibrium, player 1 of type $k$ accepts (resp.,
rejects) player 2's proposal $x$ when $U^{k}(x)>u_{0}^{k}$ (resp., $%
U^{k}(x)<u_{0}^{k}$). We further assume that player 1 accepts the proposal
when he is indifferent, which is consistent with our interest in situations
in which player 2 strictly prefers that player 1 participates. By proceeding
backwards, $\Gamma (v_{0})$ amounts to a standard sender-receiver game, with
the following utility functions (in which $I$ denotes the indicator
function):%
\begin{equation}
U_{+}^{k}(x)=U^{k}(x)I(U^{k}(x)\geq
u_{0}^{k})+u_{0}^{k}I(U^{k}(x)<u_{0}^{k})=\max \left\{
U^{k}(x),u_{0}^{k}\right\}  \label{PayU+}
\end{equation}%
for player 1 of type $k$ and%
\begin{equation}
W^{k}(v_{0},x)=V^{k}(x)I(U^{k}(x)\geq u_{0}^{k})+v_{0}I(U^{k}(x)<u_{0}^{k})%
\text{.}  \label{Payv0}
\end{equation}%
for player 2, when player 1 is of type $k$.\footnote{%
This observation is made in Chen, Nartik and Sobel (2008), in their account
of Matthews (1989).} In the latter sender-receiver game, the receiver's
utility function is not necessarily continuous, but it is
upper-semi-continuous.

\begin{lem}\label{prop1}
For every $k\in K$ and $v_{0}\in \mathbb{R}$, the utility function $%
W^{k}(v_{0},\cdot )$ defined by (\ref{Payv0}) is upper-semi-continuous.
\end{lem}

\noindent{\bf Proof:} Let $x_{n}\in X$, $x_{n}\rightarrow x$.\ The only
possibly delicate case is when $U^{k}(x_{n})<u_{0}^{k}$ for every $n$ and $%
U^{k}(x)=u_{0}^{k}$. Then $W^{k}(v_{0},x_{n})=v_{0}\leq
V^{k}(x)=W^{k}(v_{0},x)$, using our assumption. $\blacksquare $

\bigskip

Having determined player 1's behavior at the approval stage, we can define a
strategy for player 1 (in $\Gamma (v_{0})$ and $\Gamma $) as a mapping $%
\sigma :K\rightarrow \Delta (M)$.\ We interpret $\sigma (k)(m)$ as the
probability that player 1 sends message $m$ when his type is $k$, and denote
it as $\sigma (m\mid k)$. We adopt the following notations:%
\begin{equation}
\text{For every }m\in M\text{, }P_{\sigma }(m)=\sum\limits_{k}p^{k}\sigma
(m\mid k)\text{.}  \label{prob(m)}
\end{equation}%
\begin{equation}
\text{For every }k\in K\text{ and }m\in M\text{ s.t. }P_{\sigma }(m)>0\text{%
, }p_{m}^{k}(\sigma )=\frac{p^{k}\sigma (m\mid k)}{P_{\sigma }(m)}\text{.}
\label{posterior}
\end{equation}%
$p_{m}^{k}(\sigma )$ is thus the posterior probability of type $k$ computed
from $p$ and $\sigma $; let $p_{m}(\sigma )=(p_{m}^{k}(\sigma ))_{k\in K}$
denote the corresponding posterior probability distribution over $K$. We
have $\sum\limits_{m}P_{\sigma }(m)p_{m}(\sigma )=p$.

We say that $\sigma $ is nonrevealing\ if player 1 sends his message in a
type-independent way, namely, if $\sigma (m\mid k)=\sigma (m\mid k^{\prime
}) $ for every $m\in M$, $k$, $k^{\prime }\in K$.\ In this case, $%
p_{m}(\sigma )=p$ for every $m$ s.t. $P_{\sigma }(m)>0$.

For player 2, a strategy is a mapping $\tau :M\rightarrow X$, namely, a
\textquotedblleft pure\textquotedblright\ strategy with respect to the set $%
X $ (but as indicated above, $X=\Delta (A)$ for a finite set of actions $A$
is a particular case).\footnote{%
Restriction to pure strategies of player 2 is justified by the fact that
these will be enough to establish existence of equilibria in $\Gamma $.}

We say that $(\sigma ,\tau )$ is \textquotedblleft without
exit\textquotedblright\ if%
\begin{equation}
U^{k}(\tau (m))\geq u_{0}^{k}\text{ }\forall k\in K\text{, }\forall m\in M%
\text{ s.t. }\sigma (m\mid k)>0\text{,}  \label{Noexit}
\end{equation}%
namely, if%
\begin{equation*}
U_{+}^{k}(\tau (m))=U^{k}(\tau (m))\text{ }\forall k\in K\text{, }\forall
m\in M\text{ s.t. }\sigma (m\mid k)>0\text{.}
\end{equation*}%
Recalling (\ref{X(L)}) and denoting by supp $q$ the support of a probability
distribution $q\in \Delta (K)$, condition (\ref{Noexit}) is equivalent to%
\begin{equation}
\tau (m)\in X(\text{supp }p_{m}(\sigma ))\text{ }\forall m\in M\text{ s.t. }%
P_{\sigma }(m)>0\text{.}  \label{Noexitbis}
\end{equation}

At an equilibrium of $\Gamma $, we require that player 2's expected utility
be finite ($>-\infty $). Hence an equilibrium of $\Gamma $ is necessarily
without exit. By contrast, an equilibrium of $\Gamma (v_{0})$ may involve
exit of some types.

\begin{proposition}\label{prop2}
For every $v_{0}\in \mathbb{R}$, the game $\Gamma (v_{0})$ has a
nonrevealing\ equilibrium (possibly with exit). The game $\Gamma $ has a
nonrevealing equilibrium if and only if $X(K)\neq \emptyset $. Hence $\Gamma 
$ may not have any nonrevealing equilibrium.
\end{proposition}

\noindent \textbf{Proof:} The following strategies define a nonrevealing
equilibrium in $\Gamma (v_{0})$: player 1 sends the same message $m\in M$
whatever his type and then accepts $x$ if and only if $U^{k}(x)\geq
u_{0}^{k} $; whatever the message, player 2 chooses $x^{\ast }\in X$ to
maximize $\sum_{k}p^{k}W^{k}(v_{0},x)$, which is well-defined thanks to
Lemma 1.

In $\Gamma $, if $X(K)\neq \emptyset $, a nonrevealing equilibrium can be
achieved as above, provided that player 2 chooses $x^{\ast }\in X$ to
maximize $\sum_{k}p^{k}V^{k}(x)$ subject to $x\in X(K)$. If $X(K)=\emptyset $
and player 1's message is type-independent, condition (\ref{Noexitbis})
cannot be satisfied. $\blacksquare $

\bigskip

The next two propositions (which are established in Section 6.1) give us
some foundations to study the equilibria of the limit game $\Gamma $ by
making precise relationships between the latter and the equilibria without
exit of $\Gamma (v_{0})$.

\begin{proposition}\label{prop3}
Let $(\sigma ,\tau )$ be an equilibrium without exit in $\Gamma (v_{0})$,
for some $v_{0}\in \mathbb{R}$.\ Then $(\sigma ,\tau )$ is an equilibrium
without exit in $\Gamma (z_{0})$ for every $z_{0}\in \mathbb{R}$ such that $%
z_{0}\leq v_{0}$ and is also an equilibrium in $\Gamma $, with the same
interim expected utility as in $\Gamma (v_{0})$ for both players.
\end{proposition}

In other words, if $\Gamma $\ has no a equilibrium (which indeed may happen,
see Section 4.3), then, \textit{whatever }$v_{0}\in \mathbb{R}$\textit{, }$%
\Gamma (v_{0})$\textit{\ has no equilibrium without exit, that is, all
equilibria of }$\Gamma (v_{0})$\textit{\ must involve non-participation of
at least one type.}

\begin{proposition}\label{prop4}
Let $(\sigma ,\tau )$ be an equilibrium in $\Gamma $. Then there exists $%
v_{0}\in \mathbb{R}$ such that, for every $z_{0}\leq v_{0}$, $(\sigma ,\tau
) $ is an equilibrium without exit of $\Gamma (z_{0})$, with the same
interim expected utility as in $\Gamma $ for both players.
\end{proposition}

The previous properties are useful under a \textit{mechanism design}
perspective.\ Assume player 2 is a \textquotedblleft
principal\textquotedblright\ who cannot commit to a mechanism $\mu
:K\rightarrow X$ but receives a message from the agent (player 1) and then,
makes a decision in $X$ subject to the agent's participation constraints.\
With this interpretation, which turns out to be an extreme case of Bester
and Strausz (2001)'s model, an optimal mechanism amounts to an equilibrium
of $\Gamma (v_{0})$ which maximizes player 2's ex ante expected utility. At
a given $v_{0}$, the equilibrium of $\Gamma (v_{0})$ that is best for player
2 may involve the exit of some types of player 1.\ By contrast, player 2's
best equilibrium payoff in the limit game $\Gamma $, when it exists, is
achieved at an equilibrium $(\sigma ^{\ast },\tau ^{\ast })$ without exit.
Let us denote player 2's corresponding payoff as $v_{NE}^{\ast }$. By
Proposition 4, if $v_{0}$ is small enough, $(\sigma ^{\ast },\tau ^{\ast })$
is an equilibrium \textit{without exit} in $\Gamma (v_{0})$, giving the same
expected utility $v_{NE}^{\ast }$ to player 2. By Proposition 3, $%
v_{NE}^{\ast }$ is the \textit{best} equilibrium payoff player 2 can achieve
at an equilibrium \textit{without exit} in $\Gamma (v_{0})$. Moreover, as we
show in details in Section 6.1, if $v_{0}$ is sufficiently low, player 2
cannot expect an expected utility higher than $v_{NE}^{\ast }$ at an
equilibrium of $\Gamma (v_{0})$ \textit{with exit} of some types.\footnote{%
The observation that the principal's ex ante expected utility is maximized
when all types of the agent participate is also made in Bester and Strausz
(2001), footnote 8.} Summing up, if the limit game $\Gamma $ has an
equilibrium, an optimal mechanism, when the principal's utility $v_{0}$ in
case of exit is sufficiently low, can be found by maximizing player 2's
utility over the equilibria of $\Gamma $, without worrying about the precise
level $v_{0}$.

\section{Existence of a partitional equilibrium in $\Gamma $}

In this section, we focus on the game $\Gamma $.\ We identify various
sufficient conditions for the existence of a \textit{partitional}
equilibrium in $\Gamma $, in which player 1 uses a pure strategy, namely a
mapping $\sigma :K\rightarrow M$, to send his message to player 2. The
strategy $\sigma $ then induces the partition $\left\{ K_{m}\text{, }m\in
\sigma (K)\right\} $ of $K$, with $K_{m}=\sigma ^{-1}(m)=\left\{ k\in
K:\sigma (k)=m\right\} $.\ In this case, (\ref{prob(m)}) and (\ref{posterior}%
) become respectively:%
\begin{equation}
\text{For every }m\in M\text{, }P_{\sigma }(m)=\sum\limits_{k\in K_{m}}p^{k}%
\text{.}  \label{prob(m) pure}
\end{equation}%
\begin{equation}
\text{For every }k\in K\text{ and }m\in M\text{ s.t. }P_{\sigma }(m)>0\text{%
, }p_{m}^{k}(\sigma )=\frac{p^{k}I(k\in K_{m})}{P_{\sigma }(m)}\text{.}
\label{posterior pure}
\end{equation}

As in Section 2, for player 2, we focus on strategies of the form $\tau
:M\rightarrow X$. At an equilibrium of $\Gamma $, given player 1's strategy $%
\sigma $ and the message $m$ he receives, player 2 updates his belief over $%
K $ into $p_{m}(\sigma )$.$\ $To avoid exit, we must have

\begin{equation}
X(\text{supp }p_{m}(\sigma ))\neq \emptyset  \label{Xsupp}
\end{equation}%
Player 2's strategy $\tau $ is then a best response to $\sigma $ in $\Gamma $
if and only if%
\begin{equation}
\forall m\in M\text{ s.t. }P_{\sigma }(m)>0\text{, }\tau (m)\in \arg
\max_{x\in X(\text{supp }p_{m}(\sigma ))}\sum\limits_{k}p_{m}^{k}(\sigma
)V^{k}(x)\text{.}  \label{Const Opt}
\end{equation}%
We refer to these conditions as to constrained optimization.

Player 1's equilibrium conditions reduce to incentive compatibility
conditions expressing that given player 2's strategy $\tau $, player 1 of
type $k$ prefers to send $\sigma (k)$ than any other message $m$, namely,

\begin{equation}
U^{k}\left( \tau (\sigma (k))\right) )\geq U^{k}\left( \tau (m)\right) \text{
for every }k\in K\text{ and }m\in M\text{.}  \label{IC}
\end{equation}

\subsection{Two types}

If only two types are possible, we show that either there is a decision
giving both types at least their reservation utility or full revelation of
information is credible and allows to avoid exit.

\begin{proposition}\label{prop5}
Let us assume that $\mid K\mid =2$.\ Then $\Gamma $ has a partitional
equilibrium.
\end{proposition}

\noindent \textbf{Proof: }If\textbf{\ }$X(1)\cap X(2)\neq \emptyset $, let%
\begin{equation*}
x^{\ast }\in \arg \max_{x\in X(1)\cap X(2)}\left[ p^{1}V^{1}(x)+p^{2}V^{2}(x)%
\right]
\end{equation*}%
and let $m^{\ast }$ be an arbitrary element of $M$.\ Then $\sigma (1)=\sigma
(2)=m^{\ast }$ and $\tau (m)=x^{\ast }$ for every $m\in M$ defines a
nonrevealing equilibrium of $\Gamma $.

Otherwise, if $X(1)\cap X(2)=\emptyset $, let%
\begin{equation*}
x_{k}\in \arg \max_{x\in X(k)}V^{k}(x)\text{\quad }k=1,2
\end{equation*}%
and let $m_{1}\neq m_{2}$ be two distinct elements of $M$.\ Then $\sigma
(k)=m_{k}$, $\tau (m_{k})=x_{k}$, $k=1,2$, defines a fully revealing
equilibrium of $\Gamma $. Indeed, constrained optimization (\ref{Const Opt})
holds by construction; to see that incentive compatibility (\ref{IC}) also
holds, observe that $x_{k}\in X(k)$.\ Hence $x_{k}\notin X(\ell )$ for $\ell
\neq k$, since $X(1)\cap X(2)=\emptyset $.\ In other words, $U^{\ell
}(x_{k})<u_{0}^{\ell }\leq U^{\ell }(x_{\ell })$ for $\ell \neq k$. $%
\blacksquare $

\subsection{Straightforward partitional equilibria}

In this section, we propose an easy generalization of Proposition 5 when the
sender has an arbitrary number of types. Recalling the definition of $X(L)$
(see (\ref{X(L)})), let us set%
\begin{equation*}
\mathcal{T}=\left\{ \emptyset \neq L\subseteq K:X(L)\neq \emptyset \right\}
\end{equation*}%
and define $\mathcal{T}^{\ast }$ as the set of maximal elements of $\mathcal{%
T}$ for set inclusion, namely,%
\begin{equation}
\mathcal{T}^{\ast }=\left\{ L\in \mathcal{T}:\left[ L^{\prime }\in \mathcal{T%
}\text{ and }L\subseteq L^{\prime }\right] \Rightarrow L^{\prime }=L\right\} 
\text{.}  \label{Maximal elements}
\end{equation}%
We refer to $\mathcal{T}^{\ast }$ as to the \textquotedblleft participation
structure\textquotedblright\ of the game $\Gamma $.

\begin{proposition}\label{prop6}
If the participation structure of $\Gamma \mathcal{\ }$is a partition of $K$%
, $\Gamma $ has a partitional equilibrium.
\end{proposition}

\noindent \textbf{Proof:}

Let $\mathcal{T}^{\ast }=\left\{ K_{r}\right\} $.\ Consider the strategy of
player 1 consisting of revealing the cell $K_{r}$ containing his type. Let $%
x_{r}^{\ast }\in X(K_{r})$ be an optimal decision of player 2 when he learns
that player 1's type belongs $K_{r}$, namely,%
\begin{equation*}
x_{r}^{\ast }\in \arg \max_{x\in X(K_{r})}\sum_{k\in K_{r}}\frac{p^{k}}{%
\sum_{j\in K_{r}}p^{j}}V^{k}(x)\text{.}
\end{equation*}%
Constrained optimization (\ref{Const Opt}) holds by construction.\ Incentive
compatibility (\ref{IC}) is also immediate, because if $k\in $ $K_{r}$, $%
x_{r}^{\ast }\in X(K_{r})$ while for $j\neq r$, $x_{j}^{\ast }\notin
X(K_{r}) $. $\blacksquare $

\subsection{Decision in \textbf{a real interval and monotonic utility
function for the sender}}

The next result holds in particular when player 2's decision can be
interpreted as a probability distribution over two possible actions (i.e., $%
X=\Delta (A)$,$\mid A\mid =2$) and the utility $U^{k}(x)$ of player 1 of
type $k$ is expected utility with respect to\ $x$.

\begin{proposition}\label{prop7}
Let us assume that \textit{the decision set }$X$ \textit{is} \textit{a real
interval and every utility function }$U^{k}$, $k\in K$, \textit{is monotonic
over} $X$.\ Then $\Gamma $ has a partitional equilibrium.
\end{proposition}

\noindent \textbf{Proof:}

Let us take, without loss of generality, $X=\left[ 0,1\right] $; define%
\begin{eqnarray*}
K_{-} &=&\left\{ k\in K:U^{k}\text{ is weakly decreasing and not constant}%
\right\} \\
K_{+} &=&\left\{ k\in K:U^{k}\text{ is weakly increasing or constant}\right\}
\end{eqnarray*}%
We can find $x_{0}^{k}$, $k\in K$, such that if $k\in K_{-}$, $U^{k}(x)\geq
u_{0}^{k}\Leftrightarrow x\leq x_{0}^{k}$ and if $k\in K_{+}$, $U^{k}(x)\geq
u_{0}^{k}\Leftrightarrow x\geq x_{0}^{k}$. We define next

\begin{center}
$%
\begin{array}{cccc}
x_{-}=\min_{k\in K_{-}}x_{0}^{k} & \text{if }K_{-}\neq \emptyset \text{ ;} & 
x_{-}=1 & \text{if }K_{-}=\emptyset \text{.} \\ 
x_{+}=\max_{k\in K_{+}}x_{0}^{k} & \text{if }K_{+}\neq \emptyset \text{ ;} & 
x_{+}=0 & \text{if }K_{+}=\emptyset \text{.}%
\end{array}%
$
\end{center}

\noindent If $x_{+}\leq x_{-}$ (in particular, if $K_{-}$ or $K_{+}$ $%
=\emptyset $), let%
\begin{equation*}
x^{\ast }\in \arg \max_{\left[ x_{+},x_{-}\right] }\sum_{k\in K}p^{k}V^{k}(x)
\end{equation*}%
and let $m^{\ast }$ be an arbitrary element of $M$.\ Then $\sigma
(k)=m^{\ast }$ for every $k\in K$ and $\tau (m)=x^{\ast }$ for every $m\in M$
defines a nonrevealing equilibrium of $\Gamma $.

\noindent If $x_{+}>x_{-}$, let $m_{-}^{\ast }\neq m_{+}^{\ast }$ be two
distinct elements of $M$.\ Take $\sigma (k)=m_{-}^{\ast }$ if $k\in K_{-}$, $%
\sigma (k)=m_{+}^{\ast }$ if $k\in K_{+}$, namely, $\sigma $ induces the
partition $\left\{ K_{-},K_{+}\right\} $.\ Player 2's corresponding
posterior probability distribution on $K$ can be computed as in (\ref%
{posterior pure}):%
\begin{equation*}
p_{m_{-}^{\ast }}^{k}=\frac{p^{k}I(k\in K_{-})}{\sum_{j\in K_{-}}p^{j}}\text{
and }p_{m_{+}^{\ast }}^{k}=\frac{p^{k}I(k\in K_{+})}{\sum_{j\in K_{+}}p^{j}}%
\text{.}
\end{equation*}%
Let then%
\begin{equation*}
x_{-}^{\ast }\in \arg \max_{\left[ 0,x_{-}\right] }\sum\limits_{k\in
K_{-}}p_{m_{-}^{\ast }}^{k}V^{k}(x)\text{ and }x_{+}^{\ast }\in \arg \max_{%
\left[ x_{+},1\right] }\sum\limits_{k\in K_{+}}p_{m_{+}^{\ast }}^{k}V^{k}(x)%
\text{.}
\end{equation*}%
Constrained optimization (\ref{Const Opt}) holds by construction.\ There
remains to check incentive compatibility (\ref{IC}). Observe that $%
x_{-}^{\ast }<x_{+}^{\ast }$; for $k\in K_{-}$, $U^{k}$ is decreasing, hence 
$U^{k}(x_{-}^{\ast })\geq U^{k}(x_{+}^{\ast })$. Similarly for $k\in K_{+}$, 
$U^{k}$ is increasing so that $U^{k}(x_{-}^{\ast })\leq U^{k}(x_{+}^{\ast })$%
. $\blacksquare $

\subsection{Type-independent utility function for the receiver}

In this section, we assume that, when the sender accepts the receiver's
proposal, the utility function of the receiver does not depend on the
sender's type, namely, that $V^{k}(x)=V(x)$ for every $k$ and $x$. This
assumption is sometimes referred to as \textquotedblleft private
values\textquotedblright\ or \textquotedblleft known-own
payoff.\textquotedblright\ Matthews (1989) formulates it in the context of a
game of information transmission with sender's approval.

\begin{theorem}\label{prop8}
\textit{Let us assume that player 2's utility function does not depend on
player 1's type, namely, that there exists a continuous function }$%
V:X\rightarrow \mathbb{R}$ such that\textit{\ }$V^{k}=V$ \textit{for every} $%
k\in K$\textit{.\ }Then $\Gamma $ has a partitional equilibrium.
\end{theorem}

The proof consists of an algorithm, which constructs a partitional
equilibrium that is as revealing as possible, given the incentive
compatibility constraints to be fulfilled.\ More precisely, the initial
candidate is the fully revealing equilibrium. Imagine that type $k$ would
envy type $\ell $ if one tried to implement fully revealing strategies,
while type $\ell $ would not envy any type. By merging type $\ell $ and type 
$k$, one reduces the incentives problem. A key property is that, if player
2's utility function is independent of player 1's type, then player 2's
optimal decision $x_{\ell }$ when facing type $\ell $ remains optimal when
facing type $\ell $ \textit{or} type $k$.\ Before making use of it, we first
show, by relying on the same kind of argument, that the envy relation cannot
have any cycle.

\bigskip

\noindent \textbf{Proof:}

Let us fix, for every $k\in K$,%
\begin{equation}
x_{k}\in \arg \max_{x\in X(k)}V(x)\text{.}  \label{opti k}
\end{equation}%
The existence of such $x_{k}$'s is guaranteed by our assumptions. If the
previous optimization problem has several solutions, we take $x_{k}$ to
maximize $U^{k}(x)$.

For every pair of types\textit{\ }$j,k\in K$, we say that \textquotedblleft
type $k$ envies type $j$\textquotedblright\ -- and write $k\mathcal{R}j$ --
if $U^{k}(x_{j})>U^{k}(x_{k})$. An immediate property is that%
\begin{equation}
\text{for every }j,k\in K\text{, }k\mathcal{R}j\Rightarrow V(x_{k})>V(x_{j})%
\text{.}  \label{propR}
\end{equation}
To show this, observe that, by definition, $x_{k}\in X(k)$, i.e., $%
U^{k}(x_{k})\geq u_{0}^{k}$.\ Hence, if $k\mathcal{R}j$, we must have $%
U^{k}(x_{j})>u_{0}^{k}$, which implies $x_{j}\in X(k)$ (so that $x_{j}\in
X(j)\cap X(k)$) and $V(x_{k})\geq V(x_{j})$.\ But $V(x_{k})=V(x_{j})$ cannot
arise, because $U^{k}(x_{j})>U^{k}(x_{k})$ and, in case of multiple
solutions to $\max_{x\in X(k)}V(x)$, we choose $x_{k}$ to maximize $U^{k}(x)$%
.

The previous property implies that the envy relation $\mathcal{R}$ has no
cycle.

We will gradually construct a subset $L\subseteq K$ of \textit{leader types}
which do not envy any other type in $L$ and a subset $F=K\setminus L$ of 
\textit{follower types }which envy a type in $L$.

We start with $L=F=\emptyset $. Let us denote as $\alpha _{1}<\cdots <\alpha
_{n}$ the distinct values among $V(x_{k})$, $k\in K$. Necessarily, $n\leq
\mid K\mid $.\ Define then%
\begin{equation}
K_{j}=\left\{ k\in K:V(x_{k})=\alpha _{j}\right\} \text{\quad }j=1,...,n%
\text{.}  \label{Kj}
\end{equation}

\begin{description}
\item[Step 1] Consider every type $k\in K_{1}$: $V(x_{k})=\alpha _{1}$ is
strictly below any other $\alpha _{j}$.\ By (\ref{propR}), type $k$ cannot
envy any other type. We modify $L$ into $L=K_{1}$, while $F$ does not change
($F=\emptyset $).

\item[Step 2] Consider every type $k\in K_{2}$. If $k$ does not envy any
type, put $k$ in $L$. Otherwise, again by (\ref{propR}), $k$ can only envy a
type in $L$ (as defined at the end of step 1, namely, $K_{1}$), put $k$ in $%
F $.

\item[$\cdots $] 

\item[Step \textit{j}] Let $L$ and $F$ be the sets of leaders and followers
constructed so far. $L\cup F=K_{1}\cup \cdots \cup K_{j-1}$ so that by (\ref%
{propR}) and (\ref{Kj}), types in $L\cup F$ cannot envy types in $K_{j}$.
Consider every such type $k\in K_{j}$. If $k$ envies a type in $L$, put $k$
in $F$. Otherwise, put $k$ in $L$. $L$ and $F$ are thus updated at the end
of step \textit{j}.

\item[$\cdots $] 

\item[Step\textit{\ n}] Proceed as for step \textit{j}.\ Deduce the final
sets of leaders and followers.
\end{description}

For instance, if $\mid K\mid =3$ and $\mathcal{R}$ is fully described by $3%
\mathcal{R}2$ and $2\mathcal{R}1$, the previous construction results in $%
K_{1}=\left\{ 1\right\} $, $K_{2}=\left\{ 2\right\} $, $K_{3}=\left\{
3\right\} $, $L=\left\{ 1,3\right\} $.

Using the $x_{k}$'s defined by (\ref{opti k}) and the set $L$, we construct
an equilibrium $(\sigma ,\tau )$ of $\Gamma $. For simplicity, we rename the
messages in $M$ so that $L\subseteq M$. Player 1's strategy is such that $%
\sigma (K)=L$.\ More precisely, $\sigma :K\rightarrow L$ is defined by

\begin{center}
$%
\begin{array}{lll}
\sigma (k)= & k & \text{if }k\in L \\ 
& \arg \max_{j\in L,k\mathcal{R}j}U^{k}(x_{j}) & \text{if }k\in K\setminus L%
\text{.}%
\end{array}%
$
\end{center}

\noindent In other words, leader types announce themselves, while non leader
types report the leader type they most envy.\ Player 2's strategy is defined
by $\tau :L\rightarrow X:\tau (\ell )=x_{\ell }$, with $x_{\ell }$ defined
by (\ref{opti k}).

Incentive compatibility (\ref{IC}) follows from the fact that player 2's
strategy $\tau $ restricts his decisions to the subset $\left\{ x_{\ell
},\ell \in L\right\} $. Hence types in $L$, who cannot envy any other type
in $L$, are truthful.\ Types in $K\setminus L$ behave as well as they can
given the player 2's restricted decision set.

If player 1 follows $\sigma $, then, given message $\ell \in L$, player 2
deduces that player 1's type $k\in \sigma ^{-1}(\ell )$.\ The set $\sigma
^{-1}(\ell )$ contains $\ell $, $x_{\ell }\in X(\ell )$ by (\ref{opti k})
and all other types in $\sigma ^{-1}(\ell )$ envy $\ell $, so that $x_{\ell
}\in \bigcap_{k\in \sigma ^{-1}(\ell )}X(k)$.\ Since $x_{\ell }$ is a
maximizer of $V(x)$ over $X(\ell )$, it is also a maximizer of $V(x)$ over $%
\bigcap_{k\in \sigma ^{-1}(\ell )}X(k)$. $\blacksquare $

\bigskip

\noindent \textbf{Remarks:}

\begin{description}
\item[-] A main feature of the proof of Theorem 8 is that, in the
partitional equilibrium that is constructed, the receiver makes a decision
in a subset of $\left\{ x_{k},k\in K\right\} $ where $x_{k}$ is the optimal
decision he would make if he were sure to face type $k$. The receiver's
private values guarantee that if type $k$ envies type $\ell $, then $x_{\ell
}$, the receiver's optimal choice when he faces type $\ell $ for sure (i.e.,
under the constraint $x\in X(\ell )$), is still optimal when he faces type $%
k $ or type $\ell $ (i.e., under the constraint $x\in X(\ell )\cap X(k)$).\
This property may no longer hold when player 2's utility is type-dependent.

\item[-] Theorem 8 does not depend on the underlying utility representation:
the result holds if the receiver's von Neumann-Morgenstern \textit{%
preferences} over $X$ given type $k$ are equivalent for every $k\in K$.
\end{description}

\section{Examples (including a counter-example)}

\subsection{Partitional equilibrium}

Let the informed player have three possible types, i.e., $K=\left\{
1,2,3\right\} $ and let the uninformed player's decision set be%
\begin{equation*}
X=\left\{ (x_{a},x_{b}):x_{a}\geq 0,x_{b}\geq 0,x_{a}+x_{b}\leq 100\right\} 
\text{.}
\end{equation*}%
Let the utility function and reservation utility of the informed player be

\bigskip

\begin{center}
$%
\begin{array}[t]{cc}
U^{1}(x)=x_{a}-x_{b} & u_{0}^{1}=30\text{,} \\ 
U^{2}(x)=x_{b}-x_{a} & u_{0}^{2}=40\text{,} \\ 
U^{3}(x)=x_{a}+2x_{b} & u_{0}^{3}=20\text{.}%
\end{array}%
$

\bigskip
\end{center}

\noindent Let the uninformed player's utility function be type-independent:%
\begin{equation*}
V^{k}(x)=V(x)=-(x_{a}+x_{b})\text{, }k=1,2,3\text{.}
\end{equation*}

There are two goods, $a$ and $b$, $X$ accounts for the decision-maker's
resource constraints.\ Type 1 likes good $a$, dislikes good $b$; type 2 has
symmetric preferences; type 3 likes both goods, and likes good $b$ more than
good $a$.

As a possible interpretation, player 1 is a kidnapper who can have political
motivations (type 1), just look for a monetary ransom (type 2) or be
opportunistic (type 3).\ Good $a$ stands for political prisoners who can be
released while good $b$ stands for money.\ If player 1 does not accept
player 2's offer, the hostage is killed, leading to an invaluable loss for
player 2.

Recalling (\ref{X(L)}) and using \textquotedblleft $Co$\textquotedblright\
for convex hull, we have here%
\begin{eqnarray*}
X(\left\{ 1\right\} ) &=&X(\left\{ 1,3\right\} )=Co\left\{
(30,0),(100,0),(65,35)\right\} \text{,} \\
X(\left\{ 2\right\} ) &=&X(\left\{ 2,3\right\} )=Co\left\{
(0,40),(0,100),(30,70)\right\} \text{,} \\
X(\left\{ 3\right\} ) &=&Co\left\{ (20,0),(100,0),(0,100),(0,10)\right\} 
\text{,} \\
X(\left\{ 1,2\right\} ) &=&X(\left\{ 1,2,3\right\} )=\emptyset \text{.}
\end{eqnarray*}

Assume first that player 1's type $k$ is known, namely, that $p^{k}=1$. Let
then $x_{k}^{\ast }$ be the uninformed player's optimal decision (in $X$)
when he faces type $k$:%
\begin{equation}
x_{1}^{\ast }=(30,0)\text{, }x_{2}^{\ast }=(0,40)\text{, }x_{3}^{\ast
}=(0,10)\text{.}  \label{CR ex}
\end{equation}

Suppose next that only two types are possible.\ If $p^{1}=0$, given that $%
X(\left\{ 2,3\right\} )\neq \emptyset $, there is a nonrevealing
equilibrium, $x_{2}^{\ast }=(0,40)$.\ Similarly for $p^{2}=0$, with $%
x_{1}^{\ast }=(30,0)$.\ If $p^{3}=0$, there is no way to satisfy type 1 and
type 2 at the same time.\ But there is a completely revealing equilibrium: $%
x_{1}^{\ast }=(30,0)$ to type 1, $x_{2}^{\ast }=(0,40)$ to type 2 is
incentive compatible. This illustrates Proposition 5.

Let the three types be possible, namely $p^{k}>0$ for $k=1,2,3$. There is no
nonrevealing equilibrium, since $X(\left\{ 1,2,3\right\} )=\emptyset $.\
There is no completely revealing equilibrium either: (\ref{CR ex}) implies
that type 3 would pretend to be type 2 (type 3 envies type 1 and type 2 even
more).

As expected from Theorem 8, there exists a partitional equilibrium.\ The
informed player is invited to report whether his type is 1 or not. If he
reports type 1, the uninformed player proposes $x_{1}^{\ast }=(30,0)$.\ If
the informed player reports that his type is not 1, the uninformed player
proposes%
\begin{equation*}
\arg \min_{x\in X(\left\{ 2,3\right\} )}(x_{a}+x_{b})=x_{2}^{\ast }=(0,40)%
\text{.}
\end{equation*}%
As in the proof of Theorem 8, incentive compatibility is ensured by the fact
that the decision proposed to type 3 is the one he most envies among $%
x_{1}^{\ast }$ and $x_{2}^{\ast }$.

\subsection{Mixed equilibrium}

Let us modify the uninformed player's utility function in the previous
example, to make it depend on the informed player's type:%
\begin{eqnarray*}
V^{1}(x) &=&\frac{x_{a}}{3}\text{,} \\
V^{2}(x) &=&\frac{x_{b}}{3}\text{,} \\
V^{3}(x) &=&-(x_{a}+x_{b})\text{.}
\end{eqnarray*}%
A possible interpretation is that the decision-maker is happy to pay\ when
the kidnapper has sharp preferences.

Let as above $x_{k}^{\ast }$ denote the uninformed player's optimal decision
(in $X$) when he faces type $k$; we have now%
\begin{equation}
x_{1}^{\ast }=(100,0)\text{, }x_{2}^{\ast }=(0,100)\text{, }x_{3}^{\ast
}=(0,10)\text{.}  \label{CR exbis}
\end{equation}

Let us take $p=(\frac{1}{3},\frac{1}{3},\frac{1}{3})$.\ There is no\textit{\ 
}nonrevealing equilibrium, since $X(\left\{ 1,2,3\right\} )=\emptyset $.
There is no\textit{\ }completely revealing equilibrium: given (\ref{CR exbis}%
), type 3 would pretend to be type 2. More generally, there is no
partitional equilibrium.\ Given the above description of the sets $X(L)$,
two possible partitions must still be considered: $\left\{ \left\{ 1\right\}
,\left\{ 2,3\right\} \right\} $ and $\left\{ \left\{ 1,3\right\} ,\left\{
2\right\} \right\} $.

\bigskip

\noindent $\left\{ \left\{ 1\right\} ,\left\{ 2,3\right\} \right\} $: if the
uninformed player believes he faces type 1 (posterior $(1,0,0)$), his
optimal choice is $x_{1}^{\ast }=(100,0)$; if he believes he faces type 2 or
type 3 (posterior $(0,\frac{1}{2},\frac{1}{2})$), his optimal choice is 
\begin{equation*}
x_{23}^{\ast }=\arg \min_{x\in X(\left\{ 2,3\right\} )}\left[ x_{a}+\frac{2}{%
3}x_{b}\right] =(0,40)\text{.}
\end{equation*}%
This cannot be incentive compatible for type 3:%
\begin{equation*}
100=U^{3}(x_{1}^{\ast })>U^{3}(x_{23}^{\ast })=80\text{.}
\end{equation*}

\bigskip

\noindent $\left\{ \left\{ 1,3\right\} ,\left\{ 2\right\} \right\} $: if the
uninformed player believes he faces type 2 (posterior $(0,1,0)$), his
optimal choice is $x_{2}^{\ast }=(0,100)$; if he believes he faces type 1 or
type 3 (posterior $(\frac{1}{2},0,\frac{1}{2})$), his optimal choice is 
\begin{equation*}
x_{13}^{\ast }=\arg \min_{x\in X(\left\{ 1,3\right\} )}\left[ \frac{2}{3}%
x_{a}+x_{b}\right] =(100,0)\text{.}
\end{equation*}%
Again, this cannot be incentive compatible for type 3:%
\begin{equation*}
200=U^{3}(x_{2}^{\ast })>U^{3}(x_{13}^{\ast })=100\text{.}
\end{equation*}%
This illustrates that \textit{Theorem 8 does not extend to the case where
player 2's utility function depends on player 1's type.}

Let us show that if player 1 uses a \textit{mixed} strategy, a partially
revealing equilibrium exists in this example: type 1 reports that his type
belongs to $\left\{ 1,3\right\} $, type 2 reports that his type belongs to $%
\left\{ 2,3\right\} $, type 3 reports that his type belongs to $\left\{
1,3\right\} $ (resp., $\left\{ 2,3\right\} $) with probability $\frac{1}{3}$
(resp., $\frac{2}{3}$). If the informed player follows this reporting
strategy, the uninformed player's posterior upon receiving $\left\{
1,3\right\} $ is $(\frac{3}{4},0,\frac{1}{4})$ while upon receiving $\left\{
2,3\right\} $, it is $(0,\frac{3}{5},\frac{2}{5})$. Given $\left\{
1,3\right\} $, the uninformed player's problem reduces to $\min_{x\in
X(\left\{ 1,3\right\} )}x_{b}$.\ Every $x=(x_{a},0)$ with $x_{a}\in \left[
30,100\right] $ is optimal.\ Let us take $x_{13}^{\ast }=(80,0)$. Given $%
\left\{ 2,3\right\} $, the uninformed player's optimal choice is%
\begin{equation*}
x_{23}^{\ast }=\arg \min_{x\in X(\left\{ 2,3\right\} )}\left[ 2x_{a}+x_{b}%
\right] =(0,40)\text{.}
\end{equation*}%
There remains to check incentive compatibility.\ Type 1 prefers $%
x_{13}^{\ast }=(80,0)$ to $x_{23}^{\ast }=(0,40)$, and vice-versa for type
2. Type 3 must be indifferent between sending $\left\{ 1,3\right\} $ or $%
\left\{ 2,3\right\} $, because he must randomize between these two
outcomes.\ Indeed we have $U^{3}(x_{13}^{\ast })=U^{3}(x_{23}^{\ast })=80$.
Proposition 9 in Section 5 states that the previous construction can be
generalized.

\subsection{No equilibrium at all}

In the following example, none of the existence results of Section 3 can be
applied. We will show that there is no equilibrium, even if player 1 makes
use of a mixed strategy.$\ $The game is described by:$\mid K\mid =3$, $%
X=\Delta (A)$, where $A=\left\{ a,b,c\right\} $, $u_{0}^{k}=0$, $k=1,2,3$.\
The following payoff matrices describe $(U^{k}(\alpha ),V^{k}(\alpha ))$ for
every $\alpha \in A$:

\bigskip

\begin{center}
$%
\begin{array}{cccc}
& a & b & c \\ 
&  &  &  \\ 
k=1 & 0,2 & -2,0 & 1,1 \\ 
&  &  &  \\ 
k=2 & 1,1 & 0,2 & -2,0 \\ 
&  &  &  \\ 
k=3 & -2,0 & 1,1 & 0,2%
\end{array}%
$

\bigskip
\end{center}

The utility functions over $X=\Delta (A)$ are obtained as expected utilities
with respect to $x=(x_{a},x_{b},x_{c})$.

If player 2 knows that he faces type $k$ (i.e., $p^{k}=1$), he gets his
first best by choosing $a$ if $k=1$, $b$ if $k=2$, $c$ if $k=3$. But if $%
p^{k}>0$ for every $k$, there is no nonrevealing equilibrium ($%
\sum_{k}U^{k}(x)<0$ for every $x\in \Delta (A)$) and no fully revealing
equilibrium (incentive compatibility is violated).

Looking for a partially revealing equilibrium, we first check that there is
a unique, nonrevealing equilibrium, as soon as only two types are possible.\
Take, e.g., $p^{3}=0$. Then%
\begin{eqnarray*}
X(\left\{ 1,2\right\} ) &=&\left\{ x\in X:-2x_{b}+x_{c}\geq 0\text{ and }%
x_{a}-2x_{c}\geq 0\right\} \\
&=&Co\left\{ (1,0,0),(\frac{2}{3},0,\frac{1}{3}),(\frac{4}{7},\frac{1}{7},%
\frac{2}{7})\right\}
\end{eqnarray*}%
and player 2's optimization problem is:%
\begin{equation*}
\max p^{1}(2x_{a}+x_{c})+p^{2}(x_{a}+2x_{b})\text{ s.t. }x\in X(\left\{
1,2\right\} )\text{.}
\end{equation*}%
For every $p$ such that $p^{1}>0$ and $p^{2}>0$, the unique solution is
achieved at $x=(1,0,0)$, namely, action $a$ with probability 1. Similarly,
action $b$ is the only solution if $p^{2}>0$ and $p^{3}>0$, action $c$ is
the only solution if $p^{1}>0$ and $p^{3}>0$.

Let us start with $p$ such that $p^{k}>0$ for every $k$.\ By sending his
message according to a mixed strategy $\sigma $, player 1 \textquotedblleft
splits\textquotedblright\ the prior belief $p$ into posteriors $p_{m}(\sigma
)$ such that $\sum\limits_{m}P_{\sigma }(m)p_{m}(\sigma )=p$ (see (\ref%
{posterior})).\ Taking account of incentive compatibility, $p$ cannot be
split (only) extreme points, because there is no fully revealing
equilibrium.\ At least one of the posteriors $p_{m}$ must be on an edge, say 
$p_{m}^{3}=0$, so that $\tau (m)=a$. There should be another posterior $%
p_{m^{\prime }}$ with $p_{m^{\prime }}^{3}>0$, with $\tau (m^{\prime })=b$
or $c$. To achieve the posteriors $p_{m}$ and $p_{m^{\prime }}$, message $m$
must be sent with positive probability by types 1 and 2, while message $%
m^{\prime }$ must be sent with positive probability by at least type 3. If $%
\tau (m^{\prime })=b$, type 2 strictly prefers $m^{\prime }$ to $m$.\ If $%
x(m^{\prime })=c$, type 1 strictly prefers $m^{\prime }$ to $m$. Hence there
is no incentive compatible splitting and thus no equilibrium at all, even if
player 1 can use a mixed strategy.

A mediated equilibrium, in which information transmission is monitored by a
mediator, can nevertheless be achieved in the previous example.\ Consider
the following three lotteries over $A$: $\delta ^{1}=(\frac{1}{2},0,\frac{1}{%
2})$, $\delta ^{2}=(\frac{1}{2},\frac{1}{2},0)$, $\delta ^{3}=(0,\frac{1}{2},%
\frac{1}{2})$.\ Assume that, instead of selecting a message by himself,
player 1 can just choose among these three lotteries. If player 1 expects
player 2 to pick the action selected by the lottery, player 1 prefers $%
\delta ^{k}$ over the other two lotteries when his type is $k$.\ Similarly,
player 2 is happy to choose the action recommended by the lottery if he
believes that player 1 reveals his type truthfully to the mediator. This
procedure will be generalized in the next section.

\section{Equilibrium in the case of three types}

In this section, we propose a thorough analysis of the equilibria of $\Gamma 
$ when player 1 has three possible types\textit{\ }($\mid K\mid =3$). Recall
that $\mathcal{T}^{\ast }$ denotes the participation structure of $\Gamma $
(see (\ref{X(L)}) and (\ref{Maximal elements})).

\bigskip

When $\mid K\mid =3$, there are three typical cases:

\begin{enumerate}
\item $\mathcal{T}^{\ast }$ is a partition of $K$.

\item $\mathcal{T}^{\ast }=\left\{ \left\{ 1,2\right\} ,\left\{ 1,3\right\}
,\left\{ 2,3\right\} \right\} $.

\item $\mathcal{T}^{\ast }=\left\{ \left\{ 1,3\right\} ,\left\{ 2,3\right\}
\right\} $.
\end{enumerate}

In case 1, by Proposition 6, $\Gamma $ has\textit{\ }a partitional
equilibrium. Case 2 means that player 2 is able to obtain the approval of
every pair of types but cannot ensure the participation of the three types
simultaneously.\ In this case, as illustrated in Section 4.3, $\Gamma $ may
have no mixed equilibrium.\ We will show below that $\Gamma $ always has a 
\textit{mediated} equilibrium. In case 3, which has been illustrated in
Section 4.2., player 2 can only guarantee the approval of two pairs of
types.\ We will establish that under further assumptions, $\Gamma $ always
has then a mixed (possibly not pure) equilibrium.

\subsection{$\mathcal{T}^{\ast }=\left\{ \left\{ 1,2\right\} ,\left\{
1,3\right\} ,\left\{ 2,3\right\} \right\} $}

Let us enrich the description of $\Gamma $ by adding a mediator who invites
player 1 to report a type (in $K$) and then selects a decision (in $X$) that
he recommends player 2.\ At a mediated equilibrium, player 1 truthfully
reveals his type to the mediator, player 2 proposes to player 1 the decision 
$x$ that is recommended by the mediator and finally, player 1 of type $k$
accepts player 2's proposal $x^{\prime }$ provided that $U^{k}(x^{\prime
})\geq u_{0}^{k}$.

Let us construct a mediated equilibrium in $\Gamma $ when $K=\left\{
1,2,3\right\} $ and $\mathcal{T}^{\ast }=\left\{ \left\{ 1,2\right\}
,\left\{ 1,3\right\} ,\left\{ 2,3\right\} \right\} $.\ Recall that $p^{k}>0$
for every $k$. For every pair $(j,k)$ of types, let $x_{jk}^{\ast }\in
X(\left\{ j,k\right\} )$ be an optimal decision for player 2 when he learns
that player 1 is of type $j$ or $k$, namely,%
\begin{equation}
x_{jk}^{\ast }\in \arg \max \left[ \frac{p^{j}}{p^{j}+p^{k}}V^{j}(x)+\frac{%
p^{k}}{p^{j}+p^{k}}V^{k}(x)\right] \text{.}  \label{Const opt ex}
\end{equation}%
Consider the following mediator: for every $k=1,2,3$, if player 1 reports
type $k$, he selects $x_{ik}^{\ast }$ or $x_{jk}^{\ast }$, $i,j\neq k$, $%
i\neq j$, with equal probability $\frac{1}{2}$ and recommends it to player
2.\ If player 1 reports his type $k$ truthfully, player 2 learns, with equal
probability, that player 1's type is in $\left\{ i,k\right\} $ or in $%
\left\{ j,k\right\} $ for $i,j\neq k$, $i\neq j$.\ Condition (\ref{Const opt
ex}) guarantees that player 2 follows the mediator's recommendation.

For player 1, let us consider $k=1$.\ By reporting his type truthfully,
player 1 obtains%
\begin{equation*}
\frac{1}{2}U^{1}(x_{12}^{\ast })+\frac{1}{2}U^{1}(x_{13}^{\ast })\text{.}
\end{equation*}%
If he lies by, say, pretending to be of type $2$, he obtains\footnote{%
Player 1's incentive compatibility condition reflects the fact that this
player can veto player 2's offer, namely, can lie about his type and/or
reject player 2's proposal. This is a \textquotedblleft veto-incentive
compatibility condition\textquotedblright\ (see, e.g., Forges (1999)), which
implies posterior individual rationality.}%
\begin{equation*}
U_{2}^{1}=\frac{1}{2}\max \left\{ U^{1}(x_{12}^{\ast }),u_{0}^{1}\right\} +%
\frac{1}{2}\max \left\{ U^{1}(x_{23}^{\ast }),u_{0}^{1}\right\} \text{.}
\end{equation*}%
By construction, $x_{12}^{\ast }\in X(\left\{ 1,2\right\} )$ and $%
x_{13}^{\ast }\in X(\left\{ 1,3\right\} )$.\ Hence $U^{1}(x_{12}^{\ast
})\geq u_{0}^{1}$ and $U^{1}(x_{13}^{\ast })\geq u_{0}^{1}$. But $%
x_{23}^{\ast }\notin X(1)$ because $\left\{ 2,3\right\} $ is maximal. Hence, 
$U^{1}(x_{23}^{\ast })<u_{0}^{1}$.\ 
\begin{equation*}
U_{2}^{1}=\frac{1}{2}U^{1}(x_{12}^{\ast })+\frac{1}{2}u_{0}^{1}\leq \frac{1}{%
2}U^{1}(x_{12}^{\ast })+\frac{1}{2}U^{1}(x_{13}^{\ast })\text{.}
\end{equation*}%
The other incentive compatibility conditions of player 1 can be checked in a
symmetric way. $\blacksquare $

\subsection{$\mathcal{T}^{\ast }=\left\{ \left\{ 1,3\right\} ,\left\{
2,3\right\} \right\} $}

\begin{proposition}\label{prop9}
Let us assume that $K=\left\{ 1,2,3\right\} $, the participation structure $%
\mathcal{T}^{\ast }=\left\{ \left\{ 1,3\right\} ,\left\{ 2,3\right\}
\right\} $, \textit{the decision set }$X$ \textit{is compact and convex} and 
\textit{the utility functions }$U^{k}$ and $V^{k}$, $k\in K$, \textit{are
affine}.\ Then $\Gamma $ \textit{has a partially revealing equilibrium}, in
which player 1 uses a mixed, possibly not pure, strategy.
\end{proposition}

\noindent{\bf Proof:} See Section 6.2.\ We establish that there must exist an
equilibrium in which type 1 reports that his type belongs to $\left\{
1,3\right\} $, type 2 reports that his type belongs to $\left\{ 2,3\right\} $
and type 3 reports that his type belongs to $\left\{ 1,3\right\} $ (resp., $%
\left\{ 2,3\right\} $) with some probability $\delta \in (0,1)$ (resp., $%
1-\delta $). Incentive compatibility requires that type 3 be indifferent
between reporting $\left\{ 1,3\right\} $ or $\left\{ 2,3\right\} $. Such an
equilibrium is shown to be the only possible one in the example of Section
4.2.

\section{Appendix}

\subsection{Proof of Propositions 3 and 4}

For the sake of completeness, we first explicitly recall the conditions to
be satisfied by an equilibrium of $\Gamma (v_{0})$ whether they involve exit
of some types on equilibrium path or not.

Let us fix a pair of strategies $\sigma :K\rightarrow \Delta (M)$ and $\tau
:M\rightarrow X$. Player 1's equilibrium conditions\textbf{\ }can be written
as%
\begin{equation}
U_{+}^{k}(\tau (m))\geq U_{+}^{k}(\tau (m^{\prime }))\text{ }\forall k\in K%
\text{, }\forall m\in M\text{ s.t. }\sigma (m\mid k)>0,\forall m^{\prime
}\in M\text{.}  \label{EqP1}
\end{equation}

\noindent Player 2's equilibrium conditions can be written as%
\begin{equation*}
\sum\limits_{k}p_{m}^{k}(\sigma )W^{k}(v_{0},\tau (m))\geq
\sum\limits_{k}p_{m}^{k}(\sigma )W^{k}(v_{0},x)\text{ }\forall m\in M\text{
s.t. }P_{\sigma }(m)>0\text{, }\forall x\in X\text{.}
\end{equation*}

We deduce that\textbf{\ }the necessary and sufficient conditions for $%
(\sigma ,\tau )$ to be an equilibrium without exit are:

\bigskip 

For player 1:%
\begin{equation}
U^{k}(\tau (m))\geq U^{k}(\tau (m^{\prime }))\text{ }\forall k\in K\text{, }%
\forall m\in M\text{ s.t. }\sigma (m\mid k)>0,\forall m^{\prime }\in M\text{,%
}  \label{EqP1noexit}
\end{equation}
implying that%
\begin{equation*}
U^{k}(\tau (m))=U^{k}(\tau (m^{\prime }))\text{ }\forall m,m^{\prime }\in M%
\text{ s.t. }\sigma (m\mid k)>0\text{ and }\sigma (m^{\prime }\mid k)>0\text{%
.}
\end{equation*}

For player 2:%
\begin{equation}
\tau (m)\in \left[ \arg \max_{x\in X}\sum\limits_{k}p_{m}^{k}(\sigma
)W^{k}(v_{0},x)\right] \cap X(\text{supp }p_{m}(\sigma ))\text{,}
\label{OptiNoexit}
\end{equation}%
implying constrained optimization, namely (\ref{Const Opt}).

\bigskip

Recall that by definition, an equilibrium of $\Gamma $ cannot involve exit,
so that the conditions for $(\sigma ,\tau )$ to be an equilibrium in $\Gamma 
$ are thus just (\ref{EqP1noexit}) and (\ref{Const Opt}).

\bigskip

\noindent \textbf{Remarks on} (\ref{EqP1noexit}):

\begin{itemize}
\item As a refinement of subgame perfect equilibrium, Matthews (1989)
strengthens (\ref{EqP1}) to (\ref{EqP1noexit}) in the case of equilibria
which typically involve exit on path.

\item For an equilibrium without exit, player 1's equilibrium conditions
take the simple form (\ref{EqP1noexit}) because player 2's strategy $\tau $
is pure.\ In an equilibrium without exit, for every $k\in K$, player 1's
equilibrium strategy consists of sending $m$, selected with probability $%
\sigma (m\mid k)$, and then accept player 2's proposal, namely, $\tau (m)$,
which is fully anticipated at the time to choose $m$.\ A deviation consists
of sending $m^{\prime }\in M$ (possibly such that $\sigma (m^{\prime }\mid
k)=0$) and then, accept or reject player 2's proposal $\tau (m^{\prime })$.\
Player 1's equilibrium conditions thus take the form%
\begin{equation*}
U^{k}(\tau (m))\geq \max \left\{ U^{k}(\tau (m^{\prime })),u_{0}^{k}\right\} 
\text{ }\forall k\in K\text{, }\forall m\in M\text{ s.t. }\sigma (m\mid
k)>0,\forall m^{\prime }\in M\text{.}
\end{equation*}%
These are equivalent to (\ref{EqP1noexit}), since (\ref{Noexit}) holds at an
equilibrium without exit. If player 2's strategy $\tau $ were mixed ($\tau
:M\rightarrow \Delta (X)$), we would have to write\footnote{%
The expression is similar to \textquotedblleft veto-incentive
compatibility\textquotedblright\ (see, e.g., Forges (1999)) and implies
posterior individual rationality.}%
\begin{equation*}
\sum\limits_{x}\tau (x\mid m)U^{k}(x)\geq \sum\limits_{x}\tau (x\mid
m^{\prime })\max \left\{ U^{k}(x)),u_{0}^{k}\right\} \text{.}
\end{equation*}
\end{itemize}

\bigskip

\noindent \textbf{Proposition 3.} \textit{Let }$(\sigma ,\tau )$\textit{\ be
an equilibrium without exit in }$\Gamma (v_{0})$\textit{, for some }$%
v_{0}\in \mathbb{R}$\textit{.\ Then }$(\sigma ,\tau )$\textit{\ is an
equilibrium without exit in }$\Gamma (z_{0})$\textit{\ for every }$z_{0}\in
R $\textit{\ such that }$z_{0}\leq v_{0}$\textit{\ and is also an
equilibrium in }$\Gamma $\textit{, with the same interim expected utility as
in }$\Gamma (v_{0})$\textit{\ for both players.}

\bigskip

\noindent \textbf{Proof:} If $(\sigma ,\tau )$ satisfies (\ref{OptiNoexit})
in $\Gamma (v_{0})$, the same holds in $\Gamma (z_{0})$ since for every $%
x\in X$,%
\begin{equation*}
\sum\limits_{k}p_{m}^{k}(\sigma )W^{k}(z_{0},x)\leq
\sum\limits_{k}p_{m}^{k}(\sigma )W^{k}(v_{0},x)
\end{equation*}%
and for every $x\in X($supp $p_{m}(\sigma ))$,%
\begin{equation*}
\sum\limits_{k}p_{m}^{k}(\sigma
)W^{k}(z_{0},x)=\sum\limits_{k}p_{m}^{k}(\sigma
)W^{k}(v_{0},x)=\sum\limits_{k}p_{m}^{k}(\sigma )V^{k}(x)\text{.}
\end{equation*}

Furthermore, (\ref{Const Opt}) must hold and player 1's equilibrium
conditions (\ref{EqP1noexit}) are the same in $\Gamma (v_{0})$ and $\Gamma
(z_{0})$ or $\Gamma $ as long as player 2's strategy remains unchanged. $%
\blacksquare $

\bigskip

\noindent \textbf{Proposition 4.} \textit{Let }$(\sigma ,\tau )$\textit{\ be
an equilibrium in }$\Gamma $\textit{. Then there exists }$v_{0}\in \mathbb{R}
$\textit{\ such that, for every }$z_{0}\leq v_{0}$\textit{, }$(\sigma ,\tau
) $\textit{\ is an equilibrium without exit of }$\Gamma (z_{0})$\textit{,
with the same interim expected utility as in }$\Gamma $\textit{\ for both
players.}

\bigskip

\noindent \textbf{Proof:} Let $(\sigma ,\tau )$ be an equilibrium without
exit in $\Gamma $.\ By definition, constrained optimization (\ref{Const Opt}%
) holds, so that in particular $\tau (m)\in X($supp $p_{m}(\sigma ))$ for
every $m$ such that $P_{\sigma }(m)>0$. Let us keep player 1's strategy, $%
\sigma $, fixed. Player 2's strategy $\tau $ remains a best reply to $\sigma 
$ in $\Gamma (v_{0})$, with $v_{0}\leq \min_{k\in K}\min_{x\in X}V^{k}(x)$,
provided that $v_{0}$ is such that optimality of no exit holds, namely,
recalling (\ref{OptiNoexit}),%
\begin{equation*}
\forall m\in M\text{ s.t. }P_{\sigma }(m)>0,\sum\limits_{k}p_{m}^{k}(\sigma
)V^{k}(\tau (m))\geq \max_{x\in X}\sum\limits_{k}p_{m}^{k}(\sigma
)W^{k}(v_{0},x)\text{.}
\end{equation*}%
These can be viewed as finitely many inequalities over $v_{0}$,which have a
solution in $\mathbb{R}$, since the RHS are well-defined, for every $%
v_{0}\leq \min_{k\in K}\min_{x\in X}V^{k}(x)$, by Lemma 1.\footnote{%
The RHS of the inequalities can be rewritten as%
\begin{equation*}
\max_{L\varsubsetneq \text{supp}(p_{m}(\sigma ))}\max_{x\in X(L)}\left\{
\sum\limits_{k\in L}p_{m}^{k}(\sigma )V^{k}(x)+v_{0}\sum\limits_{k\in \text{%
supp}(p_{m}(\sigma ))\diagdown L}p_{m}^{k}(\sigma )\right\} \text{.}
\end{equation*}%
} Hence there exists $v_{0}$ such that $(\sigma ,\tau )$ is an equilibrium
without exit of $\Gamma (v_{0})$ and from Proposition 3, in $\Gamma (z_{0})$%
, for every $z_{0}\leq v_{0}$. $\blacksquare $

\bigskip

\noindent \textbf{Application to a mechanism design problem}

In Section 3, we have established that, under various reasonable
assumptions, $\Gamma $ has a partitional equilibrium $(\sigma ,\tau )$, in
which both $\sigma $ and $\tau $ are pure.\ In this case, we can easily
compute the highest ex ante expected utility that player 2, interpreted here
as the principal, can obtain at a partitional equilibrium of $\Gamma $.%
\footnote{%
The number of pure equilibrium payoffs is finite, in the same way as the
number of partitions of $K$.\ Hence as soon as there is a pure equilibrium
in $\Gamma $, there is an equilibrium achieving the highest expected payoff
for the receiver.} As explained in Section 2, this application is inspired
by a particular case of Bester and Strausz (2001)'s model.

Let us show that there exists $v_{0}\in \mathbb{R}$ such that, for every $%
z_{0}\leq v_{0}$, the highest ex ante expected utility player 2 can obtain
at an arbitrary partitional equilibrium of $\Gamma (z_{0})$ (which can
involve exit or not) is the same as in $\Gamma $.\ More precisely, there
exists $v_{0}\in \mathbb{R}$ such that, for every $z_{0}\leq v_{0}$, the
best partitional equilibrium for player 2 in $\Gamma $ remains the best
partitional equilibrium for player 2 in $\Gamma (z_{0})$.

Let $v_{NE}^{\ast }$ be the highest ex ante expected utility player 2 can
obtain at a partitional equilibrium of $\Gamma $.\ This number is
well-defined if $\Gamma $ has a partitional equilibrium.\ Let $(\sigma
^{\ast },\tau ^{\ast })$ achieve the expected utility $v_{NE}^{\ast }$ for
player 2. Using Proposition 4, there exists $v_{0}$ sufficiently small such
that for every $z_{0}\leq v_{0}$, $(\sigma ^{\ast },\tau ^{\ast })$ is an
equilibrium without exit of $\Gamma (z_{0})$ with the same expected utility $%
v_{NE}^{\ast }$ for player 2. By Proposition 3, for every such $z_{0}$,
there does not exist any equilibrium \textit{without exit} giving a higher
expected utility to player 2 (because such an equilibrium would still be an
equilibrium of $\Gamma $, with the same expected utilities).

Let us consider the partitional equilibria $(\sigma ,\tau )$ of $\Gamma
(v_{0})$ in which exit possibly occurs, i.e., in which the set%
\begin{equation*}
K_{E}=\left\{ k\in K:U^{k}(\tau \circ \sigma (k))<u_{0}^{k}\right\} \neq
\emptyset \text{,}
\end{equation*}%
i.e., $p_{E}=_{def}\sum_{k\in K_{E}}p^{k}>0$. The highest expected utility
player 2 can achieve at such an equilibrium is%
\begin{equation*}
p_{E}v_{0}+(1-p_{E})\overline{v}
\end{equation*}%
where%
\begin{equation*}
\overline{v}=\max_{k\in K}\max_{x\in X}V^{k}(x)\text{.}
\end{equation*}%
If $v_{0}$ is such that, for every $p_{E}$ that can arise given the prior $p$%
,%
\begin{equation*}
p_{E}v_{0}+(1-p_{E})\overline{v}\leq v_{NE}^{\ast }\text{,}
\end{equation*}%
\begin{equation}
\text{namely, }v_{0}\leq \frac{1}{p_{E}}\left[ v_{NE}^{\ast }-(1-p_{E})%
\overline{v}\right]  \label{highest ut}
\end{equation}%
then $(\sigma ^{\ast },\tau ^{\ast })$ will guarantee the highest possible
equilibrium utility to player 2, in every game $\Gamma (z_{0})$ with $%
z_{0}\leq v_{0}$.

Let $\underline{k}$ be the type with the smallest prior probability, namely,
such that $p^{\underline{k}}=\min \left\{ p^{1},\cdots ,p^{K}\right\} $. The
inequality (\ref{highest ut}) will hold at every $p_{E}$ that can arise
given the prior $p$ as soon as it holds at $p_{E}=p^{\underline{k}}$: we
just have to require%
\begin{equation*}
v_{0}\leq \frac{1}{p^{\underline{k}}}\left[ v_{NE}^{\ast }-(1-p^{\underline{k%
}})\overline{v}\right] \text{.}
\end{equation*}

The previous result is quite intuitive: an upper bound on the receiver's
expected utility at an equilibrium of $\Gamma (v_{0})$ with exit is obtained
when the receiver's proposal is rejected by only the least likely type,
while the best possible utility is achieved at all the other types. If $%
v_{0} $ is sufficiently low, the best equilibrium utility for the receiver
in $\Gamma (v_{0})$ will be not be achieved at an equilibrium with exit, but
rather at an equilibrium without exit, which is in turn is necessarily an
equilibrium of $\Gamma $.

\subsection{Proof of Proposition 9}

\noindent \textbf{Proposition 9.} \textit{Let us assume that} $K=\{1,2,3\},$ ${\cal T}^*=\{\{1,3\},\{2,3\}\}$ \textit{the decision set} $X$ \textit{is compact and convex and the utility functions} $U^k$ \textit{and} $V^k, k \in K,$ \textit{are affine. Then}  $\Gamma$ \textit{has a partially revealing} {\it equilibrium}. 

\smallskip

For simplicity we assume here that $u_0^k=0$ for each $k$. This is w.l.o.g. since we can translate the payoffs of each type of the sender. We start with preliminaries. 

\subsubsection{Mappings and multi-valued mappings}

\noindent Define for each $p$ in $\Delta(K)$: 

$f(p)=\sup\{\sum_{k \in K} p^k V^k(x), x \in X(\supp p)\}\in \R\cup\{-\infty\}$, 

$Y(p)=\argmax_{x \in X(\supp p)} \sum_k p^kV^k(x)\subset X(\supp p)$, and 

$\Phi(p)=\{(U^k(x))_{k \in K}, x \in Y(p)\} \; \subset \R^K.$\\

The sets $Y(p)$ and $\Phi(p)$ are convex compact subsets of $\R$ and $\R^K$, respectively. 
If  $X(\supp p)\neq \emptyset$, then $f(p)\in \R$, $Y(p)\neq \emptyset$ and $\Phi(p)\neq \emptyset$. For each  $u\in \Phi(p)$, we have $u^k\geq 0$ for each $k\in \supp p$.
At an equilibrium of $\Gamma$, if the belief of the receiver (after having received the message of the sender) is $p$, then he has to  propose a decision in $Y(p)$, inducing a vector payoff in $\Phi(p)$ for the different types of player 1.

We will use  in the sequel the following three lemmas (Lemma \ref{lem4} is a simple mean-value theorem for correspondences). 
 
\begin{lem} \label{lem2} The mapping $f$ is u.s.c. and convex.

 If $p_n\xrightarrow[n \to \infty]{} p \in \Delta(K)$ with $\supp p_n=\supp p$ for each $n$, then $f(p_n)\xrightarrow[n \to \infty]{} f(p)$.
\end{lem}

\noindent{\bf Proof:} Suppose  $p_n\xrightarrow[n \to \infty]{} p$. Then for $n$ large enough, $\supp p_n\supset \supp p$ so $X(\supp p_n)\subset X(\supp p)$. It follows that $\limsup_n f(p_n)\leq f(p)$. (whether $f(p)=-\infty$ or not)

If $\supp p_n=\supp p$ for each $n$, then $|f(p_n)-f(p)|\leq \sup_{x \in X} \sum_{k \in K} |p_n^k-p^k| |V^k(x)|$, and $f(p_n)\xrightarrow[n \to \infty]{} f(p)$.

If $p=\lambda p_1+(1-\lambda) p_2$ with $\lambda\in (0,1)$, then $\supp p=\supp p_1 \cup  \; \supp p_2$ and  $X(\supp p)= X(\supp p_1)\cap X(\supp p_2)$. If $f(p)=-\infty$ then $\lambda f(p_1)+(1-\lambda) f(p_2)\geq f(p)$. Consider  $x$ in $X(\supp p)$, we have $f(p_1)\geq \sum_{k \in K} p_1^k V^k(x)$ and $f(p_2)\geq \sum_{k \in K} p_2^k V^k(x)$, so $\lambda f(p_1)+(1-\lambda) f(p_2)\geq \sum_{k \in K} p^k V^k(x)$,  and taking the supremum for $x$ in $X(\supp p)$ we get $\lambda f(p_1)+(1-\lambda) f(p_2)\geq f(p)$.   Hence $f$ is convex. $\blacksquare $

\begin{lem} \label{lem3} Consider a converging sequence  $p_n\xrightarrow[n \to \infty]{} p \in \Delta(K)$. 

a) Assume $\limsup_n f(p_n)=f(p)$. Then if $u_n\xrightarrow[n \to \infty]{} u \in \R^K$, with $u_n\in \Phi(p_n)$ for each $n$, we have
$u\in \Phi(p)$,

b) Otherwise $\limsup_n f(p_n)<f(p)$. Then  there exists $n_0$ such that for each $u\in \Phi(p)$  and $n \geq n_0$, one can find  $k\in \supp p_n\backslash \{\supp p\}$ such that $u^k<0$.
\end{lem}

\noindent{\bf Proof:} a) Without loss of generality  we assume that  $f(p_n)\xrightarrow[n \to \infty]{} f(p)$.  Write $u_n=(U^k(x_n))_{k \in K}$ with $x_n$ in $Y(p_n)$ for each $n$. By taking a converging subsequence we can assume that $x_n$ converges to some $x$ in $X$. Since $p_n\xrightarrow[n \to \infty]{} p $, for $n$ large enough $\supp p_n\supset \supp p$ so $x\in X(\supp p)$. And $\sum_{k \in K}p_n^k V^k(x_n)=f(p_n)\xrightarrow[n \to \infty]{} f(p)$, so $\sum_{k \in K}p^k V^k(x)=f(p)$. Then  $x$ belongs to $Y(p)$, and $u \in \Phi(p)$.

b) Assume that $\limsup_n f(p_n)<f(p)$. We first claim that for $n$ large enough, $Y(p)\cap X(\supp (p_n))=\emptyset$. Otherwise, we can find $x$ in $Y(p)\cap X(\supp p_n))$ for infinitely many $n$'s, we have $f(p)=\sum_k p^k V^k(x)$ and $f(p_n)\geq \sum_k p_n^k V^k(x)$ for infinitely many $n$'s, so $\limsup_n f(p_n)\geq f(p)$ which is a contradiction. We have shown that there exists $n_0$ such that for $n\geq n_0$, $Y(p)\cap X(\supp (p_n))=\emptyset$. If $x\in Y(p)$ and $n\geq n_0$, then $x\notin X(\supp p_n)$. So if $u\in \Phi(p)$ and $n\geq n_0$, there exists $k\in \supp p_n\backslash \{\supp p\}$ such that $u^k<0$.$\blacksquare $

\begin{lem} \label{lem4} Let $F:[0,1] \rightrightarrows \R$ be a correspondence with non empty convex values and compact graph. If $F(0)\subset \{x \in \R, x< 0\}$ and $F(1)\subset \{x \in \R, x> 0\}$, there exists $t$ in $(0,1)$ such that $0\in F(t)$.\end{lem}

\noindent{\bf Proof:}  The sets $C_+=\{t\in [0,1], F(t)\cap \R_+\neq \emptyset\}$ and $C_-=\{t\in [0,1], F(t)\cap \R_-\neq \emptyset\}$ are closed because $F$ is u.s.c. Since $F$ has non empty  values, $C_+$ and $C_-$ are   non empty, and  $C_+ \cup C_-=[0,1]$.   By connexity of $[0,1]$, one can find $t$ in both sets, that is such that $F(t)$ intersects both $\R_+$ and $\R_-$. Since $F(t)$ is convex, it contains  0. $\blacksquare $

\subsubsection{\bf Existence of an equilibrium}

For $k=2,3$, define  $\delta_k$ as the Dirac measure on the state $k$, and $p_{-k}$ as the conditional  probability on $K$ knowing the state is not $k$: 
$$\delta_2=(0,1,0),\;\;  \delta_3=(0,0,1), \;\;p_{-2}=(\frac{p_1}{p_1+p_3},0, \frac{p_3}{p_1+p_3})\;\; {\rm and} \;\; p_{-3}=(\frac{p_1}{p_1+p_2}, \frac{p_2}{p_1+p_2},0).$$
Choose $u_2$ in $\Phi(\delta_2)$, $u_3$ in $\Phi(\delta_3)$, $u_{1,2}$ in $\Phi(p_{-3})$ and $u_{1,3}$ in $\Phi(p_{-2})$. These are vectors in $\R^3$, and to simply notations we write:
$$u_2=\left(
\begin{array}{c}
a\\
+ \\
- \\
\end{array}
\right), u_3=\left(
\begin{array}{c}
b\\
- \\
+ \\
\end{array}
\right),  u_{1,2}=\left(
\begin{array}{c}
c\geq 0\\
+\\
- \\
\end{array}
\right), u_{1,3}=\left(
\begin{array}{c}
d\geq 0\\
-\\
+\\
\end{array}
\right).$$
\noindent with  $a=u_2^1$, $b=u_3^1$, $c=u_{1,2}^1$ and $d=u_{1,3}^1$. 
Here $+$ means $\geq 0$, and $-$ means $<0$. We have  $u_2^2\geq 0$, $u_3^3\geq 0$, $c\geq 0$, $u_{1,2}^2\geq 0$, $d\geq 0$ and $u_{1,3}^3\geq 0$ since for each $p$ and $u\in \Phi(p)$, we have $u^k\geq 0$ for each $k\in \supp p$. The subset $\{2,3\}$ is not in ${\cal T}$, this gives $u_2^3<0$, $u_3^2<0$,  $u_{1,2}^3< 0$ and  $u_{1,3}^2< 0$.

Suppose $a\leq d$. Then a simple equilibrium exists. Player 1 uses the partition $\{\{2\}, \{1,3\}\}$ to communicate: he sends the message  $m=2$ if the state is 2, and the message  $m=\{1,3\}$ if the state is 1 or 3. Player 2 proposes $x_2$ in $Y(\delta_2)$ such that $u_2=(U^k(x_2))_{k \in K}$ after receiving   $m=2$, and proposes  $x_{1,3}$ in $Y(p_{-2})$ such that $u_{1,3}=(U^k(x_{1,3}))_{k \in K}$ after receiving $m=\{1,3\}$. By definition of $Y(\delta_2)$ and $Y(p_{-2})$, player 2 is in best reply. And no type of player 1 has an incentive to deviate, so we have an equilibrium where player 1 plays pure. If we suppose $b\leq c$,  we have a similar equilibrium where player 1 uses the partition  $\{\{3\}, \{1,2\}\}$.\\

From now on, we assume  that $a>d\geq 0$ and $b>c\geq 0$. Then $a\geq 0$. Consider any    sequence $(p_n)$ converging to the Dirac measure on state 2 such that $\supp p_n =\{1,2\}$ for each $n$. By Lemma \ref{lem3} part b),    we must have $\limsup_n f(p_n)\geq f(p)$, and since $f$ is u.s.c., $\limsup_n f(p_n)= f(p)$. This being true for any such sequence, $f(p_n)\xrightarrow[n \to \infty]{} f(p)$. That is, the restriction of $f$ to the set $\{p, \supp p\subset\{1,2\}\}$ is continuous at $\delta_2$. 
 And by  Lemma \ref{lem3} part a), the restriction of $\Phi$ to the segment $[p_{-3},\delta_2]$ has a closed graph. Similarly we have $b\geq 0$, and we can prove that the restriction of $f$ to the set $\{p, \supp p\subset\{1,3\}\}$ is continuous at $\delta_3$, and 
  the restriction of $\Phi$ to the segment $[p_{-2},\delta_3]$ has a closed graph. \\

The initial probability $p$ is on the segment $[\delta_2,p_{-2}]$, and also on the segment $[\delta_3,p_{-3}]$. For $t\in [0,1]$, define $q_t=t \delta_2+(1-t)p_{-3}$ and $q'_t$ in $[p_{-2}, \delta_3]$ such that $p$ belongs to the segment $[q_t,q'_t]$. $q'_t$ is uniquely defined for each $t$, $q'_0=\delta_3$ and $q'_1=p_{-2}$. We  are going to construct an equilibrium with posteriors  $q_t$ and $q'_t$ for some appropriate $t$. We need player 1 of type 1 to be indifferent between splitting to $q_t$ and $q'_t$.
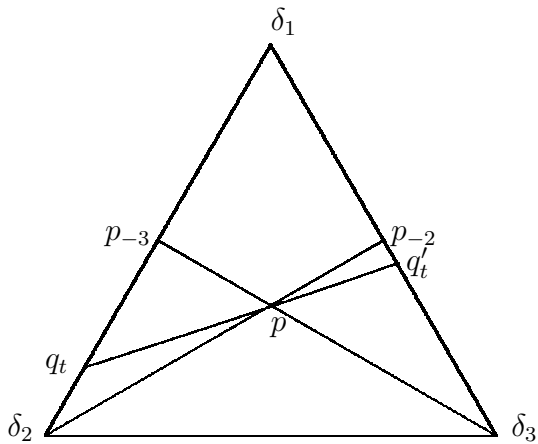
 \begin{center}
\setlength{\unitlength}{1mm}
\begin{picture}(100,60)

 \put(20,0){\line(1,0){60}}
 
 \thicklines
{ \qbezier(20,0)(35,26)(50,52)}
  {\bf  { \qbezier(80,0)(65,26)(50,52)}}
  \thinlines
  
      { \qbezier(35,26)(35,26)(80,0)}
           { \qbezier(65,26)(65,26)(20,0)}
           
            { \qbezier(25,9)(50,17)(67,23)}
   
   \put(15,0){$\delta_2$}
      \put(82,0){$\delta_3$}
        \put(50,54){$\delta_1$}
            \put(28,26){$p_{-3}$}
                  \put(66,26){$p_{-2}$}
               \put(20,9){$q_t$}         
                \put(68,22){$q'_t$}        
                
                   \put(50,14){$p$}   
\end{picture}
\end{center}

\vspace{0,5cm}

Define  the correspondence $F:[0,1] \rightrightarrows \R$, with for each $t$ in $[0,1]$:
$$F(t)=\{ u^1-v^1, u\in \Phi( q_t), v\in \Phi(q'_t)\}.$$
$F$ clearly has non empty convex compact values. We have seen that the restrictions  of $\Phi$ to the segments $[p_{-3},\delta_2]$ and $[p_{-2},\delta_3]$  have  closed graphs, moreover $q_t$ and $q'_t$ are   continuous in $t$, hence $F$ has a closed graph.
$F(0)=\{ u^1-v^1, u\in \Phi( p_{-3}), v\in \Phi(\delta_3)\}$.  If $F(0)\cap \R_+\neq \emptyset$,  there exists a  pure equilibrium where player 1 uses the partition $\{\{3\}, \{1,2\}\}$, so we assume that $F(0)$ is a subset of $\{x \in \R, x<0\}$. 
 Similarly, we assume that $F(1)=\{ u^1-v^1, u\in \Phi(\delta_2 ), v\in \Phi(p_{-2})\}$ is a subset of $\{x \in \R, x>0\}$ (otherwise there exists an equilibrium where player  1 uses the partition $\{\{2\}, \{1,3\}\}$). Then by   Lemma \ref{lem4}  we can  find $t^*$ in $[0,1]$ such that $0\in F(t^*)$.

We can now conclude the proof. We can find $x$ in $Y(q_{t^*})$, $y$ in $Y(q'_{t^*})$, $u=(U^k(x))_k\in \Phi(q_{t^*})$ and $u'=(U^k(y))_{k}\in \Phi(q'_{t^*})$ such that for some $e\geq 0$:
$$u=\left(
\begin{array}{c}
e\\
+ \\
- \\
\end{array}
\right), u'=\left(
\begin{array}{c}
e\\
- \\
+\\
\end{array}
\right)$$
We have an equilibrium as follows. Player 1 sends a message so as to induce the posteriors  $q_{t^*}$ and  $q'_{t^*}$ (type 2  sends the message 2, type 3 sends the message 3, and type 1 randomizes between the messages 2 and 3 so that the posteriors are $q_{t^*}$ after $m=2$ and  $q'_{t^*}$ after $m=3$). Player 2 then proposes $x$ at $q_{t^*}$, and $y$ at $q'_{t^*}$. Player 2 is in best reply by construction. Type 1 of player 1 is indifferent. If   type 2 of player 1  deviates and sends $q'_{t^*}$, player 2 will propose $y$ and type 2 will reject it, having the  reserve payoff of 0, which is not better than the payoff without deviating. Similarly, player 1 of type 3 has no profitable deviation, and we have an equilibrium. $\blacksquare $


\newpage

\begin{thebibliography}{99}
\bibitem{BS} Bester, H. and R. Strausz, 2001, \textquotedblleft Contracting
with imperfect commitment and the revelation principle: the single agent
case,\textquotedblright\ \textit{Econometrica} 69, 1077-1098.

\bibitem{ChenKS2008} Chen, Y., N.\ Kartik and J.\ Sobel, 2008, Selecting
cheap talk equilibria, \textit{Econometrica} 76, 117-136.

\bibitem{CompteJehiel2007} Compte, O.\ and P.\ Jehiel, 2007, On quitting
rights in mechanism design, \textit{American Economic Review Papers and
Proceedings} 97, 137-141.

\bibitem{CompteJehiel2009} Compte, O.\ and P.\ Jehiel, 2009, Veto constraint
in mechanism design: inefficiency with correlated types, \textit{American
Economic Journal: Microeconomics} 1, 182-206.

\bibitem{CS} Crawford, V. and J.\ Sobel, 1982, Strategic information
transmission, \textit{Econometrica }50, 1431-1451.

\bibitem{FF1990} Forges, F.,\ 1990, Universal mechanisms, \textit{%
Econometrica} 58, 1341-1364.

\bibitem{FF1999} Forges, F.,\ 1999, Ex post individually rational trading
mechanisms, \textit{Current trends in Economics} (edited by Alkan, A., C.
Aliprantis and N. Yannelis).

\bibitem{FH2018} Forges, F. and U.\ Horst, 2018, Sender-receiver games with
cooperation, \textit{Journal of Mathematical Economics.}

\bibitem{Gresik1991} Gresik, T., 1991, Ex ante efficient, ex post
individually rational trade, \textit{Journal of Economic Theory}.

\bibitem{MatPost1989} Matthews, S. and A. Postlewaite, 1989, Preplay
communication in two-person sealed-bid double auctions, \textit{Journal of
Economic Theory.}

\bibitem{Matthews1989} Matthews, S.,\ 1989, Veto threats: Rhetoric in a
bargaining game, \textit{Quarterly Journal of Economics }104, 347-400\textit{%
.}

\bibitem{Shimizu2013} Shimizu, T. , 2013, Cheap talk with an exit option:
the case of a discrete action space, \textit{Economics Letters }120, 397-400.

\bibitem{Shimizu2017} Shimizu, T. , 2017, Cheap talk with an exit option: a
model of exit and voice, \textit{International Journal of Game Theory }46,
1071-1088.
\end{thebibliography}
\end{document}